\documentclass[12pt]{amsart}
\usepackage{amsmath}
\usepackage{amsfonts}
\usepackage{epsfig}
\usepackage{fullpage,url}
\usepackage[dvips]{hyperref}
 \usepackage{amscd}   % for commutative diagrams
\DeclareFontEncoding{OT2}{}{} % to enable usage of cyrillic fonts

\DeclareMathOperator{\End}{End}

\DeclareMathOperator{\h}{h}
\DeclareMathOperator{\per}{per}

\DeclareMathOperator{\tor}{tor}
\DeclareMathOperator{\alg}{alg}
\DeclareMathOperator{\hhat}{\widehat{h}}

\DeclareMathOperator{\ord}{ord}

\newtheorem{theorem}{Theorem}[section]
\newtheorem{lemma}[theorem]{Lemma}

\newtheorem{proposition}[theorem]{Proposition}
\theoremstyle{definition}
\newtheorem{definition}[theorem]{Definition}

\newtheorem{conjecture}[theorem]{Conjecture}
\newtheorem{example}[theorem]{Example}

\newtheorem{Claim}[theorem]{Claim}

\theoremstyle{remark}
\newtheorem{remark}[theorem]{Remark}

\title{Points of small height on varieties defined over a function field}
\author{Dragos Ghioca}
\begin{document}

\begin{abstract}
We obtain a Bogomolov type of result for the additive group scheme in characteristic $p$. Our result is equivalent with a Bogomolov theorem for Drinfeld modules defined over a finite field.
\end{abstract}

\maketitle

\section{Introduction}
Let $A$ be an abelian variety defined over a number field $K$. We fix an algebraic closure $K^{\alg}$ for $K$ and we let $\hhat:A(K^{\alg})\rightarrow\mathbb{R}_{\ge 0}$ be the N\'{e}ron height associated to a symmetric, ample line bundle on $A$. Let $X$ be an irreducible $K^{\alg}$-subvariety of $A$. For each $n\ge 1$, we let 
\begin{equation}
\label{E:bogomolov varieties}
X_n=\left\{x\in X(K^{\alg})\mid\hhat(x)<\frac{1}{n}\right\}.
\end{equation}

The Bogomolov conjecture, which was proved in a special case by Ullmo (\cite{Ull}) and in the general case by Zhang (\cite{Zha}), asserts that if for every $n\ge 1$, $X_n$ is Zariski dense in $X$, then $X$ is the translate of an abelian subvariety of $A$ by a torsion point of $A$. A similar statement was obtained by Bombieri and Zannier (\cite{Bom}) for every power of the multiplicative group. This last paper constituted our inspiration for proving here a version of the Bogomolov conjecture for every power of the additive group in characteristic $p$.

\section{Statement of our main result}

Let $p$ be a prime number, let $\mathbb{F}_p$ be the field with $p$-elements and let $t$ be a transcendental element over $\mathbb{F}_p$. We let $\mathbb{F}_p(t)^{\alg}$ be a fixed algebraic closure of $\mathbb{F}_p(t)$ and we let $\mathbb{F}_p^{\alg}$ be the algebraic closure of $\mathbb{F}_p$ in $\mathbb{F}_p(t)^{\alg}$. We will call \emph{constants} the elements of $\mathbb{F}_p^{\alg}$.

We will construct next the usual Weil height on $\mathbb{F}_p(t)^{\alg}$.
First we construct a set of valuations on $\mathbb{F}_p(t)$ and then we extend our construction to $\mathbb{F}_p(t)^{\alg}$.

Let $R=\mathbb{F}_p[t]$. For each irreducible
polynomial $P\in R$ we let $v_P$ be the valuation on
$\mathbb{F}_p(t)$ given by 
$$v_P(\frac{Q_1}{Q_2})=\ord_P(Q_1)-\ord_P(Q_2)\text{
for every nonzero $Q_1,Q_2\in R$,}$$
where by $\ord_P(Q)$ we denote the order of the
polynomial $Q\in R$ at $P$.

Also, we construct the valuation $v_{\infty}$ on $\mathbb{F}_p(t)$
given by 
$$v_{\infty}(\frac{Q_1}{Q_2})=\deg (Q_2)-\deg
(Q_1)\text{ for every nonzero $Q_1,Q_2\in R$,}$$
where by $\deg(Q)$ we denote the degree of the
polynomial $Q\in R$.

We let $M_{\mathbb{F}_p(t)}$ be the set of all valuations $v_P$
for irreducible polynomials $P\in R$ plus the
valuation $v_{\infty}$. We let the \emph{degree} of $v_P$ be $d(v_P)=\deg(P)$ for
every irreducible polynomial $P\in R$ and we also let
$d(v_{\infty})=1$. Then, for every nonzero $x\in \mathbb{F}_p(t)$,
we have the sum formula
$$\sum_{v\in M_{\mathbb{F}_p(t)}}d(v)\cdot v(x)=0.$$

\begin{remark}
\label{R:geometrie}
The set of valuations $M_{\mathbb{F}_p(t)}$ consists of all the valuations associated to irreducible divisors of $\mathbb{P}^1_{\mathbb{F}_p(t)}$.
\end{remark}

Let $K$ be a finite extension of $\mathbb{F}_p(t)$. We let
$M_{K}$ be the set of all valuations on $K$ that
extend the valuations from $M_{\mathbb{F}_p(t)}$. We normalize
each valuation $w$ from $M_{K}$ so that the range
of $w$ is the entire $\mathbb{Z}$. For $w\in M_K$, if $v\in M_{\mathbb{F}_p(t)}$ lies below $w$, then $e(w|v)$ represents the corresponding ramification index, while $f(w|v)$ represents the relative residue degree. Also, we define 
\begin{equation}
\label{E:coherent degree}
d(w)=\frac{f(w|v)d(v)}{[K:\mathbb{F}_p(t)]}.
\end{equation}
It is easy to see that for every finite extension $K'$ of $K$ and for every $w'\in M_{K'}$, lying above $w\in M_K$, the following is true
\begin{equation}
\label{E:coherent degree general}
d(w')=\frac{f(w'|w)d(w)}{[K':K]}.
\end{equation}

Let $K$ be a finite extension of $\mathbb{F}_p(t)$ and let $x\in K$. For each $w\in M_K$, we define the function $\tilde{w}:K\rightarrow\mathbb{R}_{\le 0}$ as $\tilde{w}=\min\{w,0\}$. We define the local Weil height of $x$ at $w$ as
\begin{equation}
\label{E:local Weil height}
\h_w(x)=-d(w)\tilde{w}(x).
\end{equation}
Finally, we define the (global) Weil height of $x$ as
\begin{equation}
\label{E:global Weil height}
\h(x)=\sum_{w\in M_K}\h_w(x).
\end{equation}
The sum from \eqref{E:global Weil height} is finite because there are finitely many valuations $w\in M_K$ such that $w(x)< 0$. Also, the definition of the global Weil height of $x$ does not depend on the choice of the field $K$ containing $x$. To prove the last statement, it suffices to show the following:

\begin{equation}
\label{E:tower}
\sum_{w\in M_K}\h_w(x)=\sum_{v\in M_{\mathbb{F}_p(t,x)}}\h_v(x).
\end{equation}

For proving \eqref{E:tower}, it suffices to show that for every $v\in M_{\mathbb{F}_p(t,x)}$,

\begin{equation}
\label{E:local tower}
\sum_{\substack{w\in M_K\\w|v}}\h_w(x)=\h_v(x).
\end{equation}

By our convention on normalizing the valuations, for each $w\in M_K$, if $w|v$, then $w(x)=e(w|v)v(x)$. Thus equation \eqref{E:local tower} is equivalent with

\begin{equation}
\sum_{\substack{w\in M_K\\w|v}}d(w)e(w|v)=d(v).
\end{equation}
By the definition of $d(w)$ (see \eqref{E:coherent degree general}), the last equation is equivalent with the well-known formula 
$$\sum_{\substack{w\in M_K\\w|v}}f(w|v)e(w|v)=[K:\mathbb{F}_p(t,x)].$$

We let $h$ be the above defined Weil height on $\mathbb{F}_p(t)^{\alg}$. Let $n$ be a positive integer and let $\mathbb{G}_a^n$ be the $n$-dimensional additive group scheme. For any point $P=(x_1,\dots,x_n)\in\mathbb{G}_a^n(\mathbb{F}_p(t)^{\alg})$, for every finite field extension $K$ of $\mathbb{F}_p(t)(x_1,\dots,x_n)$ and for every $w\in M_K$, we define the local height of $P$ at $w$ as
\begin{equation}
\label{E:deflocalhhat}
\h_w(P)=\h_w(x_1,\dots,x_n)=\sum_{i=1}^n\h_w(x_i).
\end{equation}
We define the (global) height of $P$ as
\begin{equation}
\label{E:defglobalhhat}
\h(P)=\sum_{w\in M_K}\h_w(P).
\end{equation}
For the same reason as in the case of the Weil height, the height of $P$ is finite and well-defined.

\begin{theorem}
\label{T:main result}
Let $X$ be a $\mathbb{F}_p(t)^{\alg}$-subvariety of $\mathbb{G}_a^n$. Let $Y$ be the Zariski closure of the set $X(\mathbb{F}_p^{\alg})$, i.e. $Y$ is the largest $\mathbb{F}_p^{\alg}$-subvariety of $X$.

There exists a positive constant $C$, depending only on $X$, such that if $\h(P)<C$, then $P\in Y(\mathbb{F}_p(t)^{\alg})$.
\end{theorem}

We will prove Theorem \ref{T:main result} in the next section through a series of lemmas. Using our result we will prove in Section ~\ref{S:Bogomolov} a Bogomolov type theorem for Drinfeld modules defined over a finite field.

\section{Proof of our main result}

The following proposition contains standard results on the Weil height.
\begin{proposition}
\label{P:hhat properties}
For every $n\ge 1$ and every $P,Q\in\mathbb{G}_a^n(\mathbb{F}_p(t)^{\alg})$, the following statements are true:

\begin{enumerate}
\item[(i)]
$\h(P)=0$ if and only if $P\in\mathbb{G}_a^n(\mathbb{F}_p^{\alg})$.

\item[(ii)]
$\h(P+Q)\le\h(P)+\h(Q)$ (triangle inequality).

\item[(iii)]
For each $d>0$, there are only finitely many points $P$ of bounded height such that $[\mathbb{F}_p(t)(P) :\mathbb{F}_p(t)]\le d$ (Northcott theorem).
\end{enumerate}
\end{proposition}

\begin{definition}
\label{D:reduced polynomials}
We call \emph{reduced} a non-constant polynomial $f\in\mathbb{F}_{p}[t][X_1,\dots,X_n]$, whose coefficients $a_i$ have no nonconstant common divisor in $\mathbb{F}_p[t]$. If $f$ is reduced, we define its \emph{local height} at a place $v\in M_{\mathbb{F}_p(t)}$ as
$$\h_v(f)=\max_{i}\h_v(a_i).$$
Then we define the (global) height of $f$ as
$$\h(f)=\sum_{v\in M_{\mathbb{F}_p(t)}}\h_v(f).$$
\end{definition}

\begin{lemma}
\label{L:L1}
Let $f\in\mathbb{F}_{p}(t)[X_1,\dots,X_n]$ be a reduced polynomial of degree $d$. For every $k$ such that $p^k\ge 2\h(f)$, if $(x_1,\dots,x_n)\in\mathbb{G}_a^n(\mathbb{F}_p(t)^{\alg})$ satisfies $f(x_1,\dots,x_n)=0$, then either
$$\h(x_1,\dots,x_n)\ge\frac{1}{2d}$$
or
$$f(x_1^{p^k},\dots,x_n^{p^k})=0.$$
\end{lemma}

\begin{proof}
Let $k$ satisfy the inequality from the statement of Lemma \ref{L:L1}. Let $(x_1,\dots,x_n)\in\mathbb{G}_a^n(\mathbb{F}_p(t)^{\alg})$ be a zero of $f$. We let 
$$f=\sum_i a_iM_i$$
where the $a_i$'s are the nonzero coefficients of $f$ and the $M_i$'s are the monomials of $f$. For each $i$, we let $m_i=M_i(x_1,\dots,x_n)$.

Assume $f(x_1^{p^k},\dots,x_n^{p^k})\ne 0$. 

We let $K=\mathbb{F}_{p}(t)(x_1,\dots,x_n)$. If $\zeta= f(x_1^{p^k},\dots,x_n^{p^k})$, then (because $\zeta\ne 0$)
\begin{equation}
\label{E:sum formula for zeta}
\sum_{w\in M_K}d(w)w(\zeta)=0.
\end{equation}

Because $f(x_1,\dots,x_n)=0$, we get $\zeta=\zeta-f(x_1,\dots,x_n)^{p^k}$ and so,
\begin{equation}
\label{E:trick}
\zeta=\sum_i (a_i-a_i^{p^k})m_i^{p^k}.
\end{equation}

\begin{Claim}
\label{C:divizibilitate}
For every $g\in\mathbb{F}_p[t]$, $\left(t^{p^k}-t\right)\mid\left(g^{p^k}-g\right)$.
\end{Claim}

\begin{proof}[Proof of Claim \ref{C:divizibilitate}.]
Let $g=\sum_{j=0}^m b_jt^j$ with $b_j\in\mathbb{F}_p$. Then $g^{p^k}=\sum_{j=0}^m b_jt^{p^kj}$. Thus, if we prove 
\begin{equation}
\label{E:divi}
\left(t^{p^k}-t\right)\mid\left(t^{p^kj}-t^j\right)
\end{equation}
then we get the result of Claim \ref{C:divizibilitate}. 

The divisibility from \eqref{E:divi} is obvious when $j=0$ and so, we assume from now on that $j>0$. Then $t\mid\left(t^{p^kj}-t^j\right)$. Thus, because $t$ and $t^{p^k-1}-1$ are coprime in $\mathbb{F}_p[t]$, it suffices to show 
\begin{equation}
\label{E:divi2}
\left(t^{p^k-1}-1\right)\mid\left(t^{p^kj}-t^j\right).
\end{equation}
But $t^{p^kj}-t^j=t^j\left(t^{(p^k-1)j}-1\right)$. Because $\left(t^{p^k-1}-1\right)\mid\left(t^{(p^k-1)j}-1\right)$, divisibility \eqref{E:divi2} holds and so does Claim ~\ref{C:divizibilitate}.
\end{proof}

Using the result of Claim ~\ref{C:divizibilitate} and equation \eqref{E:trick}, we get
\begin{equation}
\label{E:alternative zeta}
\zeta=(t^{p^k}-t)\sum_i b_im_i^{p^k},
\end{equation}
where $b_i=\frac{a_i-a_i^{p^k}}{t^{p^k}-t}\in\mathbb{F}_p[t]$. Let $S$ be the set of valuations $w\in M_K$ such that $w$ lies above an irreducible factor (in $\mathbb{F}_p[t]$) of $t^{p^k}-t$. For each $w\in S$,

\begin{equation}
\label{E:first inequality}
d(w)\cdot w(\zeta)\ge d(w)\cdot w(t^{p^k}-t)-dp^k\h_w(x_1,\dots,x_n).
\end{equation}

For each $w\in M_K\setminus S$, because $\zeta=\sum_i a_im_i^{p^k}$,

\begin{equation}
\label{E:second inequality}
d(w)\cdot w(\zeta)\ge -\h_w(f)-dp^k\h_w(x_1,\dots,x_n).
\end{equation}

Adding all inequalities from \eqref{E:first inequality} and \eqref{E:second inequality} we obtain
\begin{equation}
\label{E:final inequality}
0=\sum_{w\in M_K}d(w)\cdot w(\zeta)\ge -\h(f)-dp^k\h(x_1,\dots,x_n)+\sum_{\substack{w\in M_K\\w(t^{p^k}-t)>0}}d(w)\cdot w(t^{p^k}-t).
\end{equation}
Just as in \eqref{E:tower} and \eqref{E:local tower}, $$\sum_{\substack{w\in M_K\\w(t^{p^k}-t)>0}}d(w)\cdot w(t^{p^k}-t))=\sum_{\substack{v\in M_{\mathbb{F}_p(t)}\\v(t^{p^k}-t)>0}}d(v)\cdot v(t^{p^k}-t)=-v_{\infty}(t^{p^k}-t)=p^k.$$
Thus, inequality \eqref{E:final inequality} yields
$$0\ge-\h(f)-dp^k\h(x_1,\dots,x_n)+p^k$$
and so,
\begin{equation}
\label{E:rough inequality}
\h(x_1,\dots,x_n)\ge\frac{1}{d}-\frac{\h(f)}{dp^k}.
\end{equation}
Because $k$ was chosen such that $p^k\ge 2\h(f)$, we conclude
\begin{equation}
\label{E:the inequality}
\h(x_1,\dots,x_n)\ge\frac{1}{2d}.
\end{equation}
\end{proof}

\begin{lemma}
\label{L:L2}
Let $k$ be a positive integer. Let $K$ be a finite field extension of $\mathbb{F}_p(t)$ and let $f\in K[X_1,\dots,X_n]$ be an irreducible polynomial. Let $k$ be a positive integer. If $$f(X_1,\dots,X_n)\mid f(X_1^{p^k},\dots,X_n^{p^k}),$$ then there exists $a\in K\setminus\{0\}$ such that $af\in \mathbb{F}_{p^k}[X_1,\dots,X_n]$.
\end{lemma}

\begin{proof}
Let $c$ be a nonzero coefficient of $f$. If we replace $f$ by $c^{-1}f$, the hypothesis of Lemma ~\ref{L:L2} applies to the new polynomial. Thus we may assume, without loss of generality that $f$ has a coefficient equal to $1$ and we will prove $f\in\mathbb{F}_{p^k}[X_1,\dots,X_n]$.

Let $K^{\per}=\bigcup_{m\ge 0} K^{1/p^m}$. There exists $m\ge 0$ and there exists $g\in K^{1/p^m}[X_1,\dots,X_n]$ such that
\begin{equation}
\label{1stprop}
f=g^{p^m}
\end{equation}
and
\begin{equation}
\label{E:2ndprop}
\text{$g$ is $K^{\per}$-irreducible.}
\end{equation}
Then $f(X_1^{p^k},\dots,X_n^{p^k})=g(X_1^{p^k},\dots,X_n^{p^k})^{p^m}$. We let $h\in K^{1/p^{k+m}}[X_1,\dots,X_n]$ such that $g(X_1^{p^k},\dots,X_n^{p^k})=h(X_1,\dots,X_n)^{p^k}$. Our hypothesis yields
\begin{equation}
\label{E:divisibility}
g(X_1,\dots,X_n)^{p^m}\mid h(X_1,\dots,X_n)^{p^{k+m}}.
\end{equation}
Using \eqref{E:2ndprop} in \eqref{E:divisibility}, we conclude $g\mid h$. But $\deg(g)=\deg(h)$ and so, there exists $u\in K^{1/p^{k+m}}$ such that $h=ug$. 

If we write $f$ as the following sum of monomials
$$f=\sum_i a_iM_i,$$
then $g=\sum_i a_i^{1/p^m}M_i^{1/p^m}$ and $h=\sum a_i^{1/p^{k+m}}M_i^{1/p^m}$. Because $f$ has a coefficient equal to $1$ (corresponding to the monomial $M$, say), both $g$ and $h$ have the same coefficient $1$ (corresponding to same monomial, $M^{1/p^m}$). This last fact shows $u=1$. Therefore, for every $i$, 
\begin{equation}
\label{E:a coefficient over finite field}
a_i^{1/p^m}=a_i^{1/p^{k+m}}.
\end{equation}
Equation \eqref{E:a coefficient over finite field} shows $f\in\mathbb{F}_{p^k}[X_1,\dots,X_n]$, as desired.
\end{proof}

\begin{lemma}
\label{L:L3}
Let $X$ be a $\mathbb{F}_p(t)^{\alg}$-variety strictly contained in $\mathbb{G}_a^n$. There exists a positive constant $C$, depending only on the defining equations for $X$, and there exists a  $\mathbb{F}_p^{\alg}$-variety $Z$, strictly contained in $\mathbb{G}_a^n$, whose defining equations depend only on the defining equations for $X$, such that for every $P\in X(\mathbb{F}_p(t)^{\alg})$, either $P\in Z(\mathbb{F}_p(t)^{\alg})$ or $\h(P)\ge C$.
\end{lemma}

\begin{proof}
Let $K$ be the smallest field extension of $\mathbb{F}_p(t)$ such that $X$ is defined over $K$. Let $p^m$ be the inseparable degree of the extension $K/\mathbb{F}_p(t)$ ($m\ge 0$). Let 
$$X_1=\bigcup_{\sigma}X^{\sigma},$$
where $\sigma$ denotes any field morphism $K\rightarrow\mathbb{F}_p(t)^{\alg}$ over $\mathbb{F}_p(t)$. The variety $X_1$ is a $\mathbb{F}_p(t)^{1/p^m}$-variety. Also, the defining equations of $X_1$ depend only on the defining equations for $X$. Thus, if we prove Lemma ~\ref{L:L3} for $X_1$, then our result will hold also for $X\subset X_1$. Hence we assume, for sake of simplicity, that $X$ is defined over $\mathbb{F}_p(t)^{1/p^m}$.

We let $\text{Frob}$ be the usual Frobenius. The variety $X'=\text{Frob}^mX$ is a $\mathbb{F}_p(t)$-variety, whose defining equations depend only on the defining equations for $X$. Assume we proved Lemma ~\ref{L:L3} for $X'$ and let $C'$ and $Z'$ be the positive constant and the $\mathbb{F}_p^{\alg}$-variety, respectively, associated to $X'$, as in the conclusion of Lemma ~\ref{L:L3}. Let $P\in X(\mathbb{F}_p(t)^{\alg})$. Then $P'=\text{Frob}^m(P)\in X'(\mathbb{F}_p(t)^{\alg})$. Thus, either 
$$\h(P')\ge C'$$ or 
$$P'\in Z'(\mathbb{F}_p(t)^{\alg}).$$
In the former case, because $\h(P)=\frac{1}{p^m}\h(P)$, we obtain a lower bound for the height of $P$, depending only on the defining equations of $X$ (note that $m$ depends only on the defining equations for $X$). In the latter case, if we let $Z$ be the $\mathbb{F}_p^{\alg}$-subvariety of $\mathbb{G}_a^n$, obtained by extracting the $p^m$-roots of the coefficients of a set of polynomials (defined over $\mathbb{F}_p^{\alg}$) which generate the vanishing ideal for $Z'$, we get the conclusion of Lemma ~\ref{L:L3}.

Thus, from now on in this proof, we assume $X$ is a $\mathbb{F}_p(t)$-variety. We proceed by induction on $n$. 

The case $n=1$ is obvious, because any subvariety of $\mathbb{G}_a$, different from $\mathbb{G}_a$, is a finite union of points. Thus we may take $Z=X(\mathbb{F}_p^{\alg})$ and let $$C=\min\left\{1,\min_{P\in X(\mathbb{F}_p(t)^{\alg})\setminus Z(\mathbb{F}_p(t)^{\alg})}\h(P)\right\}.$$ By construction, $C>0$ (there are finitely many points in $\left(X\setminus Z\right)(\mathbb{F}_p(t)^{\alg})$ and they all have positive height because for each one of them, not all coordinates are constant; if there are no points in $X(\mathbb{F}_p(t)^{\alg})\setminus X(\mathbb{F}_p^{\alg})$, then $C=1$). 

\begin{remark}
\label{R:case n=1}
The above argument proves the case $n=1$ for Theorem \ref{T:main result}, because the variety $Z$ that we chose is a subvariety of $X$.
\end{remark}

We assume Lemma \ref{L:L3} holds for $n-1$ and we prove it for $n$ ($n\ge 2$).

We fix a set of defining polynomials for $X$ which contains polynomials $P_i$ for which $$\max\deg(P_i)$$ is minimum among all possible sets of defining polynomials for $X$. We may assume all of the polynomials we chose are reduced. If all of them have coefficients from a finite field, i.e. $\mathbb{F}_p$, then Lemma ~\ref{L:L3} holds with $Z=X$ and $C$ any positive constant.

Assume there exists a reduced polynomial $f\notin\mathbb{F}_p[X_1,\dots,X_n]$ in the fixed set of defining equations for $X$. Let $f_i$ be all the $\mathbb{F}_p(t)$-irreducible factors of $f$ and let $X_i$ be the hypersurface corresponding to $f_i=0$, for each $i$. Then 
\begin{equation}
\label{E:inclusion in hypersurfaces}
X\subset\bigcup_i X_i\text{ (finite union).}
\end{equation}
The polynomials $f_i$ depend in a precise way on $f$. Thus, inclusion \eqref{E:inclusion in hypersurfaces} shows that it suffices to prove Lemma ~\ref{L:L3} for each $X_i$. For the sake of simplifying the notation, we assume $X$ is the hypersurface $f=0$ and $f$ is $\mathbb{F}_p(t)$-irreducible. Just as before, we may assume $f$ is reduced. If $f\in\mathbb{F}_p[X_1,\dots,X_n]$, Lemma ~\ref{L:L3} holds trivially for $X$ (as it was observed before). Thus we are left with the case $f\notin\mathbb{F}_p[X_1,\dots,X_n]$ is a reduced, $\mathbb{F}_p(t)$-irreducible polynomial and $X$ is the hypersurface containing the zeros of $f$.

Let $P=(x_1,\dots,x_n)\in X(\mathbb{F}_p(t)^{\alg})$. We apply Lemma ~\ref{L:L1} to $f$ and $(x_1,\dots,x_n)$ and conclude that either 
\begin{equation}
\label{E:height inequality}
\h(P)\ge\frac{1}{2\deg(f)}
\end{equation}
or there exists $k$ depending only on $\h(f)$ such that
\begin{equation}
\label{E:alternative}
f(x_1^{p^k},\dots,x_n^{p^k})=0.
\end{equation}
If \eqref{E:height inequality} holds, then we obtained a good lower bound for the height of $P$ (depending only on the defining equations for $f$). 

Assume \eqref{E:alternative} holds. Because $f$ is a reduced polynomial, whose coefficients are not all constants, i.e. from $\mathbb{F}_p$, Lemma ~\ref{L:L2} yields that $f(X_1,\dots,X_n)$ cannot divide $f(X_1^{p^k},\dots,X_n^{p^k})$. We know $f$ has more than one monomial because it is reduced and not all of its coefficients are from $\mathbb{F}_p$. Hence, without loss of generality, we may assume $f$ has positive degree in $x_n$. Because $f$ is irreducible, the resultant $R$ for the polynomials $f(X_1,\dots,X_n)$and $f(X_1^{p^k},\dots,X_n^{p^k})$ with respect to the variable $X_n$ is nonzero. Moreover, $R$ depends only on $f$.

The nonzero polynomial $R\in\mathbb{F}_p(t)[X_1,\dots,X_{n-1}]$ vanishes on $(x_1,\dots,x_{n-1})$. Applying the induction hypothesis to the hypersurface containing the zeros of $R$ in $\mathbb{G}_a^{n-1}$, we conclude there exists a $\mathbb{F}_p^{\alg}$-variety $Z$, strictly contained in $\mathbb{G}_a^{n-1}$, whose defining equations depend only on $R$ (and so, only on $X$) and there exists a positive constant $C$, depending only on $R$ (and so, only on $X$) such that either
\begin{equation}
\label{E:1stdichotomy1}
\h(x_1,\dots,x_{n-1})\ge C
\end{equation}
or
\begin{equation}
\label{E:1stdichotomy2}
(x_1,\dots,x_{n-1})\in Z(\mathbb{F}_p(t)^{\alg}).
\end{equation}
If \eqref{E:1stdichotomy1} holds, then $\h(x_1,\dots,x_{n-1},x_n)\ge\h(x_1,\dots,x_{n-1})\ge C$ and we have a height inequality as in the conclusion of Lemma ~\ref{L:L3}. If \eqref{E:1stdichotomy2} holds, then $(x_1,\dots,x_n)\in \left(Z\times\mathbb{G}_a\right)(\mathbb{F}_p(t)^{\alg})$ and $Z\times\mathbb{G}_a$ is a $\mathbb{F}_p^{\alg}$-variety, strictly contained in $\mathbb{G}_a^n$, as desired in Lemma ~\ref{L:L3}. This proves the inductive step and concludes the proof of Lemma ~\ref{L:L3}.
\end{proof}

We are ready now to prove Theorem \ref{T:main result}.
\begin{proof}[Proof of Theorem \ref{T:main result}.]
If $X=\mathbb{G}_a^n$, the conclusion is immediate. Therefore, assume from now on in this proof that $X$ is strictly contained in $\mathbb{G}_a^n$. 

We prove Theorem \ref{T:main result} by induction on $n$. The case $n=1$ was already proved during the proof of Lemma ~\ref{L:L3} (see Remark ~\ref{R:case n=1}).

We assume Theorem \ref{T:main result} holds for $n-1$ and we will prove that it also holds for $n$ ($n\ge 2$). Let $C$ and $Z$ be as in the conclusion of Lemma ~\ref{L:L3} for $X$. Also, we remember that $Y$ is the largest $\mathbb{F}_p^{\alg}$-subvariety of $X$.

Let $P=(x_1,\dots,x_n)\in X(\mathbb{F}_p(t)^{\alg})$. If $\h(P)\ge C$, then we have a height inequality as in the conclusion of Theorem ~\ref{T:main result}. 

Assume $\h(P)<C$. Then, according to Lemma \ref{L:L3}, $P\in Z(\mathbb{F}_p(t)^{\alg})$. Let $Z_1=X\cap Z$. Because $Z$ depends only on $X$, then also $Z_1$ depends only on $X$. Moreover, the $\mathbb{F}_p(t)^{\alg}$-irreducible components of $Z_1$ and $Z$ depend only on $X$. If the $\mathbb{F}_p(t)^{\alg}$-irreducible component of $Z_1$ containing $P$ equals the $\mathbb{F}_p(t)^{\alg}$-irreducible component of $Z$ containing $P$ (which is an $\mathbb{F}_p^{\alg}$-variety, because $Z$ is an $\mathbb{F}_p^{\alg}$-variety), we conclude that $P\in Y(\mathbb{F}_p(t)^{\alg})$, as desired.

Assume the geometrically irreducible component $W_1$ of $Z_1$ containing $P$ does not equal the geometrically irreducible component $W$ of $Z$ containing $P$. Because the varieties $Z$ and $Z_1$ and all their geometrically irreducible components depend only on $X$, there is a precise number of pairs $(Z_1',Z')$ such that
\begin{equation}
\label{E:condition equation 1}
Z_1'\text{ is a geometrically irreducible component of $Z_1$}
\end{equation}
and
$$Z'\text{ is a geometrically irreducible component of $Z$}$$
and
\begin{equation}
\label{E:condition equation 2}
Z_1'\cap Z'\ne\not\emptyset.
\end{equation}
Thus, the pair of varieties $(W_1,W)$, even though depends on $P$, it belongs to a prescribed set of possible pairs of varieties, depending only on $X$. We know $W_1\subset W$ and we assumed $W_1\ne W$, i.e. 
\begin{equation}
\label{E:dimension}
\dim W_1<\dim W\text{ (because both $W$ and $W_1$ are irreducible).}
\end{equation}
According to Lemma \ref{L:L3}, $\dim Z<n$ and so, $\dim W<n$. Let $d$ be the dimension of $W$ and without loss of generality, we may assume the projection $\pi:\mathbb{G}_a^n\rightarrow\mathbb{G}_a^d$, when restricted to $W$ is generically finite-to-one (after relabelling the coordinates we can achieve this anyway). 

Let $U$ be the Zariski closure of $\pi(W_1)$. Because $W_1$ is a closed subvariety of $W$ of smaller dimension, $\dim U<d$. Because $W_1$ depends only on $X$, $U$ depends only on $X$. Because $d<n$ and $U$ is a subvariety strictly contained in $\mathbb{G}_a^d$, we may apply the inductive hypothesis to $U$. Let $U_0$ be the largest $\mathbb{F}_p^{\alg}$-subvariety of $U$. We conclude there exists a positive constant $C_1$ depending only on the variety $U$ (and so, depending only on the variety $X$) such that either 
\begin{equation}
\label{E:2nddichotomy1}
\h(x_1,\dots,x_d)\ge C_1
\end{equation}
or 
\begin{equation}
\label{E:2nddichotomy2}
(x_1,\dots,x_d)\in U_0(\mathbb{F}_p(t)^{\alg}).
\end{equation}

If \eqref{E:2nddichotomy1} holds, then $\h(x_1,\dots,x_n)\ge\h(x_1,\dots,x_d)\ge C_1$. Because $C_1$ depends only on $X$, we obtain one of the two possibilities from the conclusion of Theorem ~\ref{T:main result}.

If \eqref{E:2nddichotomy2} holds, then $(x_1,\dots,x_n)\in \left(U_0\times\mathbb{G}_a^{n-d}\right)(\mathbb{F}_p(t)^{\alg})$. The $\mathbb{F}_p^{\alg}$-variety $U_0\times\mathbb{G}_a^{n-d}$ intersects $W$ in a subvariety of smaller dimension because 
$$\dim(\pi(U_0\times\mathbb{G}_a^{n-d}))=\dim(U_0)<d=\dim(\pi(W)).$$
Let $V=\left(U_0\times\mathbb{G}_a^{n-d}\right)\cap W$. By construction, $P$ lies on $V$ and $V$ is a $\mathbb{F}_p^{\alg}$-variety (both $U_0$ and $W$ are $\mathbb{F}_p^{\alg}$-varieties) of dimension strictly smaller than $W$ and implicitly, smaller than $\dim Z$. Moreover, $V$ depends only on $X$, because both $W$ and $U_0\times\mathbb{G}_a^{n-d}$ depend only on $X$.

We repeat the above analysis with $Z$ replaced by $V$ and $W$ replaced by the geometrically irreducible component $V'$ of $V$ containing $P$. Note that $\dim V'\le\dim V<\dim W$ and $V'$ belongs to a prescribed set of varieties (the geometrically irreducible components of $V$), which depend only on $X$.  If $V'\subset X$, then $P\in V'(\mathbb{F}_p(t)^{\alg})\subset Y(\mathbb{F}_p(t)^{\alg})$, which is one of the two possibilities from Theorem ~\ref{T:main result}. If not, then arguing as before, we conclude that either the height of $P$ has a positive lower bound depending only on $X$, or $P$ lies on a $\mathbb{F}_p^{\alg}$-variety of dimension strictly smaller than $\dim V$. If the former case holds, we stop because we found a lower bound for the height of $P$, depending only on $X$. If the latter case holds, then we go on as before. We already knew $\dim Z\le n-1$ and so, in $(n-1)$ steps our process must end with either finding a lower bound for $\h(P)$, depending only on $X$, or finding a $\mathbb{F}_p^{\alg}$-subvariety of $X$ containing $P$. The latter case yields $P\in Y(\mathbb{F}_p(t)^{\alg})$, as desired. 

This concludes the proof of Theorem \ref{T:main result}.
\end{proof}

\section{A Bogomolov type theorem for Drinfeld modules of finite characteristic}
\label{S:Bogomolov}

In this section we will explain how Theorem \ref{T:main result} yields a Bogomolov theorem for Drinfeld modules defined over a finite field. In order to do this, we first define the notion of a Drinfeld module.

As before, $p$ is a 
prime number. We let $q$ be a power of $p$.  
We let $C$ be a nonsingular projective curve defined
over $\mathbb{F}_q$ and we 
fix a closed point $\infty$ on $C$. Then we define $A$ as the ring of functions
on $C$ that are regular 
everywhere except possibly at $\infty$.

We let $K$ be a finitely generated field extension of $\mathbb{F}_q$. We
fix a 
morphism $i:A\rightarrow K$. We define the operator
$\tau$ as the power of the 
usual Frobenius with the property that for every $x\in
K^{\alg}$, $\tau(x)=x^q$. 
Then we let $K\{\tau\}$ be the ring of polynomials in
$\tau$ with coefficients from $K$ (the addition is the usual addition, while the multiplication is given by the usual composition of functions).

We fix an algebraic closure of $K$, denoted
$K^{\alg}$. We denote, as before, by $\mathbb{F}_p^{\alg}$ the
algebraic closure of 
$\mathbb{F}_p$ inside $K^{\alg}$. 

A Drinfeld module is a morphism $\phi:A\rightarrow
K\{\tau\}$ for which the 
coefficient of $\tau^0$ in $\phi_a$ is $i(a)$ for
every $a\in A$, and there 
exists $a\in A$ 
such that $\phi_a\ne i(a)\tau^0$. Following the
definition 
from \cite{Goss} we call $\phi$ a Drinfeld module of
generic characteristic 
if $\ker(i)=\{0\}$ and we call $\phi$ a Drinfeld
module of finite 
characteristic if $\ker(i)\ne \{0\}$. In the latter
case, we say that the 
characteristic of $\phi$ is $\ker(i)$ (which is a
prime ideal of $A$). 

We will also require the following definitions in this section.
\begin{definition}
\label{D:fieldofdef}
For a Drinfeld module $\phi:A\rightarrow K\{\tau\}$,
its field of definition is 
the smallest subfield of $K$ containing all the
coefficients of $\phi_a$, for 
every $a\in A$.
\end{definition}

\begin{definition}
Let $\phi:A\rightarrow K\{\tau\}$ be a Drinfeld module. We call the polynomial $f\in K^{\alg}\{\tau\}$ an endomorphism of $\phi$ if for every $a\in A$, $f\phi_a=\phi_af$. The set of all endomorphisms of $\phi$ is called $\End(\phi)$.
\end{definition}

\begin{definition}
Let $\phi:A\rightarrow K\{\tau\}$ be a Drinfeld module. For each nonzero $a\in A$, we let $\phi[a]$ be the set of all $x\in K^{\alg}$ such that $\phi_a(x)=0$. We define the torsion of $\phi$ be $\phi_{\tor}=\bigcup_{a\in A\setminus\{0\}}\phi[a]$.
\end{definition}

From this point on, $K$ will always be a finitely generated field of transcendence degree $1$ over $\mathbb{F}_p$, i.e. $[K:\mathbb{F}_p(t)]<\aleph_0$. We constructed in Section $2$ the sets of valuations $M_L$ for each finite extension $L$ of $\mathbb{F}_p(t)$ (so, in particular also for $K$).

Fix a non-constant $a\in A$. Let $L$ be a finite extension of $K$ and let $w\in M_L$. Remember the definition of $\tilde{w}=\min\{w,0\}$. For each $x\in L$ define 
$$V_w(x)=\lim_{n\rightarrow\infty}\frac{\tilde{w}(\phi_{a^n}(x))}{\deg\left(\phi_{a^n}\right)}.$$
The above defined function is well-defined and has the following properties (see \cite{T}):

$1)$ if $x$ and all the coefficients of $\phi_a$ are integral at $w$, then $V_w(x)=0$.

$2)$ for all $b\in A\setminus\{0\}$,
$V_w(\phi_b(x))=\deg(\phi_b)\cdot V_w(x)$. Moreover,
we can use any non-constant $a\in A$ for the
definition of $V_w(x)$ and we will always get the same
function $V_w$.

$3)$ $V_w(x\pm y)\ge\min\{V_w(x),V_w(y)\}$.

$4)$ if $x\in\phi_{\tor}$, then $V_w(x)=0$.

We define the local height of $x$ at $w$ (with respect to $\phi$) as
\begin{equation}
\label{E:localheight2}
\hhat_{\phi,w}(x)=-d(w) V_w(x).
\end{equation}

We define the (global) height of $x$ (with respect to $\phi$) as
\begin{equation}
\label{E:globalheight2}
\hhat_{\phi}(x)=\sum_{w\in M_L}\hhat_{\phi,w}(x).
\end{equation}
The global height is well defined and does not depend on the choice of the field $L$ containing $x$ (see \cite{T}). Also, $\hhat_{\phi}(x)=0$ if and only if $x\in\phi_{\tor}$ (see \cite{T}).

Let $n\ge 1$ and let $P=(x_1,\dots,x_n)\in\mathbb{G}_a^n(K^{\alg})$. We define the (global) height of $P$ with respect to the Drinfeld module $\phi$ as $$\hhat_{\phi}(P)=\sum_{i=1}^n\hhat_{\phi}(x_i).$$ Clearly, when $n=1$, this new definition coincides with the one from \eqref{E:globalheight2} (by identifying $P$ with its unique coordinate). Also, we extend the action of $\phi$ on $\mathbb{G}_a^n$ diagonally and so, a torsion point for $\phi$ in $\mathbb{G}_a^n(K^{\alg})$ is a point whose coordinates are torsion points for $\phi$ as elements of $K^{\alg}$. Using property $2)$ for $V_w$, we conclude that for every $a\in A\setminus\{0\}$,
\begin{equation}
\label{E:tranzitie}
\hhat_{\phi}(\phi_a(P))=\deg(\phi_a)\hhat_{\phi}(P).
\end{equation}

We are now ready to state the Bogomolov conjecture for Drinfeld modules.

\begin{conjecture}
\label{C:bogophi}
Let $K$ be a function field of transcendence degree $1$ over $\mathbb{F}_p$. Let $\phi:A\rightarrow K\{\tau\}$ be a Drinfeld module. Let $n\ge 1$ and let $X$ be an irreducible $K^{\alg}$-subvariety of $\mathbb{G}_a^n$. For every $m\ge 1$, we define the set
$$X_m=\left\{P\in X(K^{\alg})\mid\hhat_{\phi}(P)<\frac{1}{m}\right\}.$$
If $X_m$ is Zariski dense in $X$ for every $m\ge 1$, then there exists $\overline{\gamma}=(\gamma_1,\dots,\gamma_n)\in\phi_{\tor}^n$ and there exists a $K^{\alg}$-subvariety $Y\subset\mathbb{G}_a^n$, which is invariant under a nonzero endomorphism $f\in\End(\phi)$ such that 
$$X=\overline{\gamma}+Y.$$
\end{conjecture}

\begin{example}
\label{E:Ex1}
Assume $q=p$ and so, let $\tau$ be the usual Frobenius. Let $\lambda\in\mathbb{F}_{p^2}\setminus\mathbb{F}_p$. Let $X\subset\mathbb{G}_a^2$ be the curve $y=x^p$ and let $\phi:\mathbb{F}_p[t]\rightarrow\mathbb{F}_p(t)\{\tau\}$ be the Drinfeld module defined by $\phi_t=\lambda\tau$. Then the set of points $P\in X(\mathbb{F}_p^{\alg})$ is dense in $X$. On the other hand, each $P\in X(\mathbb{F}_p^{\alg})$ is a torsion point for $\phi$ and so, $\h(P)=0$. Hence, for each $m\ge 1$, the set $X_m$ (defined as in Conjecture ~\ref{C:bogophi}) is dense in $X$.

The variety $X$ is invariant under every power of the Frobenius, but it is not invariant under $\phi_t$ (and neither is any translate of $X$). This shows that for Drinfeld modules defined over a finite field, we cannot require the variety $Y$ from Conjecture ~\ref{C:bogophi} be invariant under the action of $\phi$.
\end{example}

\begin{example}
\label{E:Ex2}
Assume $q=p$ and so, let $\tau$ be the usual Frobenius. Let again $\lambda\in\mathbb{F}_{p^2}\setminus\mathbb{F}_p$. Let $X\subset\mathbb{G}_a^2$ be the curve $y=\lambda x$ and let $\phi:\mathbb{F}_p[t]\rightarrow\mathbb{F}_p(t)\{\tau\}$ be the Drinfeld module defined by $\phi_t=t\tau+\tau^3$. It is easy to see that $X$ is invariant under $\phi_{t^2}$, but not under $\phi_t$ (and no translate of $X$ is invariant under $\phi_t$). 

Let $P=(1,\lambda)\in X(\mathbb{F}_p(t)^{\alg})$. We can verify easily that $1\notin\phi_{\tor}$; actually, $\hhat_{\phi}(1)=\hhat_{\phi,v_{\infty}}(1)=\frac{1}{p^3}$, where $v_{\infty}$ is the valuation on $\mathbb{F}_p(t)$ associated with the negative degree of rational functions. For each $n\ge 1$, we let $P_n\in\mathbb{G}_a^2(\mathbb{F}_p(t)^{\alg})$ satisfy $\phi_{t^{2n}}(P_n)=P$. Because the first coordinate of $P$ is not a torsion point for $\phi$, all of the points $P_n$ are distinct. Morover, because $X$ is invariant under $\phi_{t^2}$, for each $n\ge 1$, $P_n\in X(\mathbb{F}_p(t)^{\alg})$. By \eqref{E:tranzitie}, $\hhat_{\phi}(P_n)=\frac{\hhat_{\phi}(P)}{p^{6n}}$. Hence the  points $P_n$ have global heights with respect to $\phi$ converging to $0$ and they are all on the curve $X$. Therefore all the sets $X_m$, as defined in Conjecture ~\ref{C:bogophi}, are Zariski dense in $X$. As explained in the above paragraph, $X$ is invariant under an endomorphism of $\phi$ (i.e. $\phi_{t^2}$), but it is not invariant under the action of $\phi$. Thus, just as in Example ~\ref{E:Ex1}, we conclude that for Drinfeld modules of finite characteristic, we cannot hope to replace the phrase "invariant under an endomorphism of $\phi$" by "invariant under the action of $\phi$" in Conjecture ~\ref{C:bogophi}.
\end{example}

Even though the above two examples show the conclusion of Conjecture ~\ref{C:bogophi} is sharp for Drinfeld modules of finite characteristic, we believe that for Drinfeld modules of generic characteristic, the phrase "invariant under an endomorphism of $\phi$" could be replaced by "invariant under the action of $\phi$" in Conjecture ~\ref{C:bogophi}. This phenomen of having to settle for less in the case of Drinfeld modules of finite characteristic is not new. For an analysis of this situation in the context of a Mordell-Lang statement for Drinfeld modules, we refer the reader to \cite{T}. The same behaviour presented in \cite{T} for the Mordell-Lang statement occurs in the context of a Manin-Mumford statement for Drinfeld modules, where once again the case of Drinfeld modules of finite characteristic requires a weakening of the conclusion. The Manin-Mumford Conjecture (stated below) for Drinfeld modules is a weak form of our Conjecture ~\ref{C:bogophi}.

\begin{conjecture}
\label{C:manmum}
Let $K$ be a function field of transcendence degree $1$ over $\mathbb{F}_p$. Let $\phi:A\rightarrow K\{\tau\}$ be a Drinfeld module. Let $n\ge 1$ and let $X$ be an irreducible $K^{\alg}$-subvariety of $\mathbb{G}_a^n$. If the set of torsion points $P\in\mathbb{G}_a^n(K^{\alg})$ for $\phi$ are dense in $X$, then there exists a point $\overline{\gamma}\subset\phi_{\tor}^n$ and there exists a $K^{\alg}$-subvariety $Y\subset\mathbb{G}_a^n$, which is invariant under a nonzero endomorphism of $\phi$, such that $X=\overline{\gamma}+Y$.
\end{conjecture}

Clearly, if the subset of torsion points of $\phi$ inside $X$ is Zariski dense in $X$, then each of the sets $X_m$ (as defined in Conjecture ~\ref{C:bogophi}) is dense in $X$. Thus, the hypothesis of Conjecture ~\ref{C:manmum} is stronger than the hypothesis of Conjecture ~\ref{C:bogophi}, while the conclusion is the same. Our Example ~\ref{E:Ex1} shows that the phrase "invariant under an endomorphism of $\phi$" cannot be replaced by "invariant under the action of $\phi$" in Conjecture ~\ref{C:manmum}, if $\phi$ is a Drinfeld module of finite characteristic. However, Thomas Scanlon proved in \cite{TS} Conjecture ~\ref{C:manmum} for Drinfeld modules of generic characteristic and he also strenghtens its conclusion by showing that the variety $Y$ from the conclusion of Conjecture ~\ref{C:manmum} is actually invariant under the action of $\phi$ on $\mathbb{G}_a^n$. This last result, together with other observations we made during the writing of the present paper led us to believe that also the conclusion of Conjecture ~\ref{C:bogophi} could be strenghten in the case of Drinfeld modules of generic characteristic.

In this paper we will prove Conjecture ~\ref{C:bogophi} for Drinfeld modules defined over a finite field.

\begin{theorem}
\label{T:finchar}
Let $n$, $K$ and $\phi$ be as in the statement of Conjecture \ref{C:bogophi} and let $X$ be any $K^{\alg}$-subvariety of $\mathbb{G}_a^n$. If $\phi$ is defined over a finite field, then Conjecture ~\ref{C:bogophi} holds.
\end{theorem}

\begin{proof}
The next lemma shows the connection between Theorem ~\ref{T:finchar} and Theorem ~\ref{T:main result}.
\begin{lemma}
\label{L:connection}
For every $P\in\mathbb{G}_a^n(K^{\alg})$, $\hhat_{\phi}(P)=\h(P)$ (as defined in \eqref{E:defglobalhhat}).
\end{lemma}

\begin{proof}[Proof of Lemma \ref{L:connection}.]
This is proved in Lemma $5.2.6$ of \cite{T}.
\end{proof}

Let $Y$ be the largest subvariety of $X$ defined over a finite field and let $C$ be the positive constant associated to $X$ as in Theorem ~\ref{T:main result}. Let $m\ge 1$ satisfy $\frac{1}{m}\le C$. Then by Theorem ~\ref{T:main result} applied to $X$, we conclude $X_m\subset Y(K^{\alg})$. Because $X_m$ is Zariski dense in $X$ and $Y$ is Zariski closed, we conclude $X=Y$. Thus $X$ is defined over a finite field and so, it is invariant under a power of the Frobenius. Because $\phi$ is defined over a finite field, there exists a power of the Frobenius in the endomorphism ring of $\phi$. Hence there exists a power $f$ of the Frobenius such that $X$ is invariant under $f$ and $f\in\End(\phi)$. Thus Theorem ~\ref{T:finchar} holds with $\overline{\gamma}=(0,\dots,0)$.
\end{proof}


\begin{thebibliography}{9}

\bibitem{Bom}
E. Bombieri, U. Zannier, \emph{Algebraic points on subvarieties of $\mathbb{G}_m^n$}. Internat. Math. Res. Notoces 1995, no. 7, 333-347.

\bibitem{T}
D. Ghioca, \emph{The arithmetic of Drinfeld modules}. PhD thesis, University of California at Berkeley, 2005.

\bibitem{Goss}
D. Goss, \emph{Basic structures of function field
arithmetic}, Ergebnisse 
der Mathematik und ihrer Grenzgebiete (3) [Results in
Mathematics and Related
  Areas (3)], 35. Springer-Verlag, Berlin, 1996.
  
\bibitem{TS}
Thomas Scanlon, \emph{Diophantine geometry of the torsion of a Drinfeld module}.  J. Number Theory  97  (2002),  no. 1, 10-25.  

\bibitem{Ull}
E. Ullmo, \emph{Positivit\'{e} et discr\'{e}tion des points alg\'{e}briques des courbes}. (French) [Positivity and discreteness of algebraic points of curves] Ann. of Math. (2) 147 (1998), no. 1, 167-179.

\bibitem{Zha}
S. Zhang, \emph{Equidistribution of small points on abelian varieties}. Ann. of Math. (2) 147 (1998), no. 1, 159-165.

\end{thebibliography}
\end{document}